\documentclass[jnth]{apjrnl}          
\usepackage{amssymb,amsmath,latexsym} 

{\theoremstyle{plain}
  \newtheorem{theorem}{Theorem}
  \newtheorem{corollary}{Corollary}
  \newtheorem{lemma}{Lemma}
}
{\theoremstyle{definition}
  \newtheorem{definition}{Definition}
}

\def\Z{\mathbb{Z}}
\def\N{\mathbb{N}}

\def\R{\mathbb{R}}
\def\C{\mathbb{C}}
\def\s{\mathfrak{s}}

\begin{document}

\title{The Frobenius Problem, Rational Polytopes, and Fourier-Dedekind Sums\thanks{Appeared in \emph{Journal of Number Theory} {\bf 96} (2002), 1--21. \\ Parts of this work appeared in the first author's Ph.D. thesis.}} 

\author{
Matthias Beck 
\\Department of Mathematical Sciences 
\\State University of New York 
\\Binghamton, NY 13902-6000
\\E-mail: matthias@math.binghamton.edu 
\and 
Ricardo Diaz 
\\Deptartment of Mathematics 
\\The University of Northern Colorado 
\\Greeley, CO 80639 
\\E-mail: rdiaz@bentley.unco.edu 
\and 
Sinai Robins 
\\Deptartment of Mathematics 
\\Temple University 
\\Philadelphia, PA 19122 
\\E-mail: srobins@math.temple.edu}

\date{October 24, 2001} 

\maketitle

\begin{abstract}  
We study the number of lattice points in integer dilates of the rational polytope 
\[ {\cal P} = \left\{ (x_{ 1 } , \dots , x_{ n }) \in \R_{ \geq 0 }^{ n }  : \sum_{k=1}^{n} x_{k} a_{k} \leq 1 \right\} , \] 
where $a_{1} , \dots , a_{n} $ are positive integers. This polytope is closely related to the 
\emph{linear Diophantine problem of Frobenius}: given relatively prime positive integers $a_{1} , \dots , a_{n}$, 
find the largest value of $t$ (the \emph{Frobenius number}) such that 
$ m_1 a_1 + \dots + m_n a_n = t $ has no  solution in positive integers $ m_{ 1 } , \dots , m_{ n } $. 
This is equivalent to the problem of finding the largest dilate $ t {\cal P}$ such that 
the facet $ \left\{ \sum_{k=1}^{n} x_{k} a_{k} = t \right\}  $ contains no lattice point.  
We present two methods for computing the Ehrhart quasipolynomials 
$ L ( \overline{\cal P} , t ) := \# ( t {\cal P} \cap \Z^{ n } ) $ and 
$ L ( {\cal P}^{\circ} , t ) := \# ( t {\cal P}^{\circ} \cap \Z^{ n } ) $. 
Within the computations a Dedekind-like finite Fourier 
sum  appears. We obtain a reciprocity law for these sums, generalizing a theorem of Gessel. 
As a corollary of our formulas, we rederive the reciprocity law for Zagier's 
higher-dimensional Dedekind sums. 
Finally, we find bounds for the Fourier-Dedekind sums and use them 
to give new bounds for the Frobenius number. 
\end{abstract}


\begin{subject}[2000 Mathematics Subject Classification] 
11D04, 05A15, 11H06 
\end{subject}


\section{\normalsize Introduction}
Let $ a_{1} , \dots , a_{n} $ be positive integers, $\Z^{ n } \subset \R^{ n } $ be the $n$-dimensional integer lattice, and 
  \begin{equation}\label{polytope} 
    {\cal P} = \left\{ (x_{ 1 } , \dots , x_{ n }) \in \R^{ n }  : \ x_{ k } \geq 0 , \sum_{k=1}^{n} a_{k} x_{k} \leq 1 , \right\} \ ,
  \end{equation} 
a rational polytope with vertices 
  \[ (0, \dots , 0), \left( \tfrac{ 1 }{ a_{ 1 } } , 0 , \dots 0 \right) , \left( 0, \tfrac{ 1 }{ a_{ 2 } }, 0, \dots , 0 \right) , \dots , \left( 0, \dots , 0, \tfrac{ 1 }{ a_{ n } }  \right) \ . \] 
For a positive integer $t \in \N $, let $ L ( \overline{\cal P} , t ) $ be the number of lattice points in the dilated polytope $ t {\cal P} = \{ t x : x \in {\cal P} \} $.  
Denote further the relative interior of ${\cal P}$ by ${\cal P}^{\circ}$ and the number of lattice points in $t {\cal P}^{\circ}$ by $ L ( {\cal P}^{\circ} , t ) $. 
Then  $ L ( {\cal P}^{\circ} , t ) $ 
and $ L ( \overline{\cal P} , t ) $ are quasipolynomials in $t$ of degree $n$ \cite{ehrhart}, i.e. expressions 
  \[ c_{n}(t) \ t^{n} + \dots + c_{1}(t) \ t + c_{0}(t) ,  \] 
where $c_{0}, \dots , c_{n}$ are periodic functions in $t$. In fact, if the $ a_{ k } $'s are pairwise relatively prime 
then $c_{1}, \dots c_{n}$ are constants,  so only $c_{0}$ will show this periodic dependency on $t$. 

Let $A = \{ a_{1} , \dots , a_{n} \}$ be a set of relatively prime positive integers, and 
  \begin{equation}\label{N_t} 
    p_A' (t) = \# \left\{ (m_{1}, \dots, m_{n}) \in \N^{n} : \sum_{k=1}^{n} m_{k} a_{k} = t \right\} . 
  \end{equation}
The function $p_A' (t)$ can be described as the number of \emph{restricted partitions of 
$t$ with parts in $A$}, where we require that each part is used at least once. 
(We reserve the name $p_A$ for the enumeration function of those partitions which do not have this restriction.) 
Geometrically, $p_A' (t) $ enumerates the lattice points on the skewed facet 
of ${\cal P}$. Define $f(a_{1}, \dots, a_{n})$ to be the largest value of $t$ for which 
  \[ p_A' (t) = 0 \ . \]
In the 19th century, Frobenius inaugurated the study of $f(a_{1}, \dots, a_{n})$. 
For $n=2$, it is known (probably at least since Sylvester \cite{sylvester}) that $f(a_{1}, a_{2}) = a_{1} a_{2}$. 
For $n>2$, all attempts for explicit formulas have proved elusive. 
Here we focus on the study of $p_A' (t) $, and show that it has an explicit 
representation as a quasipolynomial. Through the discussion of $p_A' (t) $, 
we gain new insights into Frobenius's problem. 

Another motivation to study $p_A' (t) $ is the following trivial reduction 
formula to lower dimensions: 
  \begin{equation}\label{recurs} p_{\{a_{1}, \dots, a_{n}\}}' (t) = \sum_{m>0} p_{\{a_{1}, \dots, a_{n-1}\}}' (t - m a_{n}) \ . \end{equation} 
Here we use the convention that $ p_A' (t) = 0 $ if $t \leq 0 $. This identity 
can be easily verified by viewing $ p_A' (t) $ as
  \[ p_A' (t) = \# \left\{ (m_{1}, \dots, m_{n}) \in \N^{n} : \sum_{k=1}^{n-1} m_{k} a_{k} = t - m_{n} a_{n} \right\} . \]
Hence, precise knowledge of the values of $t$ for which $ p_A' (t) = 0 $ in 
lower dimensions sheds additional light on the Frobenius number in higher dimensions. 

The number $ p_A' (t) $ appears in the lattice point count of 
${\cal P}$. It is for this reason that we decided to focus on this particular rational polytope. 
We present two methods (Sections \ref{mattSec} and \ref{sinaiSec}) for computing the terms appearing 
in $ L ( {\cal P}^{\circ} , t ) $ and $ L ( \overline{\cal P} , t ) $. 
Both methods are refinements of concepts that were earlier introduced by the authors \cite{beck,diaz}. 
In contrast to the mostly algberaic-geometric and topological ways of computing 
$ L ( {\cal P}^{\circ} , t ) $ and $ L ( \overline{\cal P} , t ) $ \cite{barv,bv,cs,guill,kk,kp}, our methods are analytic. 
In passing, we recover the Ehrhart-Macdonald reciprocity law relating $ L ({\cal P}^{\circ} , t ) $ 
and $ L ( \overline{\cal P} , t ) $ \cite{ehrhart,macdonald}. 
Within the computations a Dedekind-like finite Fourier sum 
appears, which shares some properties with its classical siblings, discussed in Section \ref{FDSum}. 
In particular, we prove two reciprocity laws for these sums: a rederivation of the reciprocity law 
for Zagier's higher-dimensional Dedekind sums \cite{zagier}, and a new reciprocity law that generalizes 
a theorem of Gessel \cite{gessel}. 
Finally, in Section \ref{lastsec} we give bounds on these generalized Dedekind sums and 
apply our results to give new bounds for the Frobenius number. 
The literature on such bounds is vast---see, for example, \cite{shockley,davison,erdos,kannan,rodseth1,rodseth2,selmer,vitek}. 


\section{\normalsize The residue method}\label{mattSec}
In \cite{beck}, the first author used the residue theorem to count lattice points in 
a lattice polytope, that is, a polytope with integer vertices. 
Here we extend these methods to the case of \emph{rational} vertices. 

We are interested in the number of lattice points in the tetrahedron ${\cal P}$ defined by (\ref{polytope}) 
and integral dilates of it. We can interpret
  \[ L ( \overline{\cal P} , t ) = \# \left\{ (m_{ 1 } , \dots , m_{ n }) \in \Z^{ n }  : m_{ k } \geq 0 , \sum_{ k=1 }^{ n } m_{ k } a_{ k } \leq t \right\} , \]
as the Taylor coefficient of $z^{ t }$ of the function
\begin{eqnarray*} &\mbox{}& \left( 1 + z^{ a_{ 1 } } + z^{ 2 a_{ 1 } } + \dots  \right) \cdots \left( 1 + z^{ a_{ n } } + z^{ 2 a_{ n } } + \dots  \right) \left( 1 + z + z^{ 2 } + \dots \right) \\
                  &\mbox{}& \quad = \frac{ 1 }{ 1 - z^{ a_{ 1 } } } \cdots \frac{ 1 }{ 1 - z^{ a_{ n } } } \frac{ 1 }{ 1 - z } \ . \end{eqnarray*} 
Equivalently 
  \begin{equation}\label{LPbar} L ( \overline{\cal P} , t ) = \mbox{Res} \left( \frac{ z^{ - t - 1 } }{ \left( 1 - z^{ a_{ 1 } }  \right) \cdots \left( 1 - z^{ a_{ n } }  \right) \left( 1 - z \right) } , z=0 \right) . \end{equation}
If this expression counts the number of lattice points in $\overline{t \cal P}$, then the 
remaining task is to compute the other residues of  
  \[ F_{ -t } (z) := \frac{ z^{ - t - 1 } }{ \left( 1 - z^{ a_{ 1 } }  \right) \cdots \left( 1 - z^{ a_{ n } }  \right) \left( 1 - z \right) } \ , \] 
and use the residue theorem for the sphere $\C \cup \left\{ \infty \right\} $. 
$F_{ -t } $ has poles at 0 and all $a_{1}^{\text{th}}, \dots, a_{n}^{\text{th}}$ roots of unity.  It is 
particularly easy to get precise formulas if the poles at the nontrivial roots of unity are 
simple. For this reason, assume in the following that $a_{1}, \dots, a_{n}$ are \emph{pairwise
relatively prime}.  Then the residues for the $a_{1}^{\text{th}}, \dots, a_{n}^{\text{th}}$ roots
of unity  are not hard to compute: Let $\lambda^{ a_{1} } = 1 \not= \lambda $, then 
\begin{eqnarray*} &\mbox{}& \mbox{Res} \left( F_{ -t } (z) , z = \lambda \right) = \\ 
                  &\mbox{}& \quad = \frac{ \lambda^{ - t - 1 } }{ \left( 1 - \lambda^{ a_{ 2 } }  \right) \cdots \left( 1 - \lambda^{ a_{ n } } \right) \left( 1 - \lambda \right) } \ \mbox{Res} \left( \frac{ 1 }{ 1 - z^{ a_{ 1 } } } , z = \lambda \right) \\ 
                  &\mbox{}& \quad = \frac{ \lambda^{ - t - 1 } }{ \left( 1 - \lambda^{ a_{ 2 } }  \right) \cdots \left( 1 - \lambda^{ a_{ n } }  \right) \left( 1 - \lambda \right) } \ \lim_{z \to \lambda} \frac{ z - \lambda }{ 1 - z^{ a_{ 1 } } } \\ 
                  &\mbox{}& \quad = - \frac{ \lambda^{ - t } }{ a_{1} \left( 1 - \lambda^{ a_{ 2 } }  \right) \cdots \left( 1 - \lambda^{ a_{ n } }  \right) \left( 1 - \lambda \right) } \ . \end{eqnarray*} 
If we add up all the nontrivial $a_{1}^{\text{th}}$ roots of unity, we obtain 

\begin{eqnarray*} &\mbox{}& \sum_{ \lambda^{ a_{1} } = 1 \not= \lambda } \mbox{Res} \left( F_{ -t } (z) , z = \lambda \right) = \\ 
                  &\mbox{}& \quad = \frac{-1}{a_{1}} \sum_{ \lambda^{ a_{1} } = 1 \not= \lambda } \frac{ \lambda^{ - t } }{ \left( 1 - \lambda^{ a_{ 2 } } \right) \cdots \left( 1 - \lambda^{ a_{ n } }  \right) \left( 1 - \lambda \right) } \\ 
                  &\mbox{}& \quad = \frac{-1}{a_{1}} \sum_{ k=1 }^{ a_{1} - 1 } \frac{ \xi^{ - k t } }{ \left( 1 - \xi^{ k a_{ 2 } } \right) \cdots \left( 1 - \xi^{ k a_{ n } }  \right) \left( 1 - \xi^{ k } \right) } \ , \end{eqnarray*} 
where $\xi$ is a primitive $a_{1}^{\text{th}}$ root of unity. This motivates the following 

\begin{definition} Let $c_{ 1 } , \dots, c_{ n } \in \Z$ be relatively prime to $c \in \Z$, and 
$t \in \Z$. Define the {\bf Fourier-Dedekind sum} as 
  \[ \sigma_{t} \left( c_{ 1 } , \dots, c_{ n } ; c  \right) = \frac{1}{c} \sum_{ \lambda^{ c } = 1 \not= \lambda } \frac{ \lambda^{ t } }{ \left( \lambda^{ c_{ 1 } } - 1 \right) \cdots \left( \lambda^{ c_{ n } } - 1 \right) } \ . \]
\end{definition} 
Some properties of $\sigma_{t} $ are discussed in Section \ref{FDSum}. 
With this notation, we can now write
\[ \sum_{ \lambda^{ a_{1} } = 1 \not= \lambda } \mbox{Res} \left( F_{ -t } (z) , z = \lambda \right) = (-1)^{n+1} \sigma_{-t} \left(a_{ 2 } , \dots, a_{ n }, 1 ; a_{ 1 }  \right) . \]
We get similar residues for the $a_{2}^{\text{th}}, \dots, a_{n}^{\text{th}}$ roots of unity. 
Finally, note that Res($F_{-t}, z=\infty$) = 0, so that the residue theorem allows us to 
rewrite (\ref{LPbar}):
\begin{theorem}\label{closure} Let ${\cal P}$ be given by {\rm (\ref{polytope})}, with $a_{ 1 } , \dots, a_{ n } $ pairwise relatively prime. Then 
  \[ L ( \overline{\cal P} , t ) = R_{ -t } \left( a_{ 1 } , \dots, a_{ n } \right) + (-1)^{n} \sum_{ j=1 }^{ n } \sigma_{-t} ( a_{ 1 } , \dots, {\hat a_{ j } } , \dots , a_{ n } , 1 ; a_{ j } )  \] 
where $ R_{ -t } \left( a_{ 1 } , \dots, a_{ n } \right) = - \mbox{\rm Res} \left( F_{ -t } (z) , z=1 \right) $, 
and $ {\hat a_{ j } } $ means we omit the term $ a_{ j } $. 
\end{theorem} 
{\it Remarks.} 1. $R_{ -t }$ can be easily calculated via 
\begin{eqnarray*} &\mbox{}& \mbox{Res} \left( F_{ -t } (z) , z=1 \right) = \mbox{Res} \left( e^{ z }  F_{ -t } (e^{ z } ) , z=0 \right)  \\
                  &\mbox{}& \quad = \mbox{Res} \left( \frac{ e^{ -t  z } }{ \left( 1 - e^{ a_{ 1 } z } \right) \cdots \left( 1 - e^{ a_{ n } z } \right) \left( 1 - e^{ z } \right) } , z=0 \right) . \end{eqnarray*} 
To facilitate the computation in higher dimensions, one can use mathematics software such as {\tt Maple} or {\tt Mathematica}. 
It is easy to see that $ R_{ -t } \left( a_{ 1 } , \dots, a_{ n } \right) $ is a polynomial in $t$ whose coefficients 
are rational expressions in $ a_{ 1 } , \dots, a_{ n } $. The first values for $R_{ -t }$ are 
\begin{eqnarray*} &\mbox{}& R_{ -t }(a) = \frac{ t }{ a } + \frac{ 1 }{ 2a } + \frac{ 1 }{ 2 } \\
                  &\mbox{}& R_{ -t }(a,b) = \frac{ t^{ 2 }  }{ 2ab } + \frac{ t }{ 2 } \left( \frac{ 1 }{ a } + \frac{ 1 }{ b } + \frac{ 1 }{ ab }  \right) + \frac{ 1 }{ 4 } \left( 1 + \frac{ 1 }{ a } + \frac{ 1 }{ b }  \right) + \\
                  &\mbox{}& \quad + \frac{ 1 }{ 12 } \left( \frac{ a }{ b } + \frac{ b }{ a } + \frac{ 1 }{ ab }  \right) \\
                  &\mbox{}& R_{ -t }(a,b,c) = \frac{ t^{ 3 }  }{ 6abc } + \frac{ t^{ 2 }  }{ 4 } \left( \frac{ 1 }{ ab } + \frac{ 1 }{ ac } + \frac{ 1 }{ bc } + \frac{ 1 }{ abc }  \right) + \\
                  &\mbox{}& \quad + \frac{ t }{ 12 } \left( \frac{ 3 }{ a } + \frac{ 3 }{ b } + \frac{ 3 }{ c } + \frac{ 3 }{ ab } + \frac{ 3 }{ ac } + \frac{ 3 }{ bc } + \frac{ a }{ bc } + \frac{ b }{ ac } + \frac{ c }{ ab } + \frac{ 1 }{ abc }  \right) + \\
                  &\mbox{}& \quad +  \frac{ 1 }{ 24 } \left( 3 + \frac{ 3 }{ a } + \frac{ 3 }{ b } + \frac{ 3 }{ c } + \frac{ a }{ b } + \frac{ a }{ c } + \frac{ b }{ a } + \frac{ b }{ c } + \frac{ c }{ a } + \frac{ c }{ b } + \frac{ 1 }{ ab } + \frac{ 1 }{ ac } + \frac{ 1 }{ bc } + \frac{ a }{ bc } + \frac{ b }{ ac } + \frac{ c }{ ab } \right) . \end{eqnarray*} 
2. If $a_{ 1 } , \dots, a_{ n } $ are \emph{not} pairwise relatively prime, we can get similar 
formulas for $L ( \overline{\cal P} , t ) $. In this case we do not have only simple poles, 
so that the computation of the residues gets slightly more complicated. 

For the computation of $L ( {\cal P}^{\circ} , t )$ (the number of lattice points in the \emph{interior} 
of our tetrahedron ${t \cal P}$), we similarly write 
\[ L ( {\cal P}^{\circ} , t ) = \# \left\{ (m_{ 1 } , \dots , m_{ n }) \in \Z^{ n }  : m_{ k } > 0 , \sum_{ k=1 }^{ n } m_{ k } a_{ k } < t \right\} \ .\]
So now we can interpret $ L ( {\cal P}^{\circ} , t ) $ as the Taylor coefficient of $z^{ t }$ of the function
\begin{eqnarray*} &\mbox{}& \left( z^{ a_{ 1 } } + z^{ 2 a_{ 1 } } + \dots  \right) \cdots \left( z^{ a_{ n } } + z^{ 2 a_{ n } } + \dots  \right) \left( z + z^{ 2 } + \dots \right) \\
                  &\mbox{}& \quad = \frac{ z^{ a_{ 1 } } }{ 1 - z^{ a_{ 1 } } } \cdots \frac{ z^{ a_{ n } } }{ 1 - z^{ a_{ n } } } \ \frac{ z }{ 1 - z } \ , \end{eqnarray*} 
or equivalently as 
\begin{eqnarray*} &\mbox{}& \mbox{Res} \left( \frac{ z^{ a_{ 1 } } }{ 1 - z^{ a_{ 1 } } } \cdots \frac{ z^{ a_{ n } } }{ 1 - z^{ a_{ n } } } \ \frac{ z }{ 1 - z } \ z^{ - t - 1 } , z=0 \right)  \\
                  &\mbox{}& \quad = \mbox{Res} \left( \frac{ -1 }{ z^{ 2 }  } \ \frac{ 1 }{ z^{ a_{ 1 } } - 1 } \cdots \frac{ 1 }{ z^{ a_{ n } } - 1 } \ \frac{ 1 }{ z-1 } \ z^{ t + 1 }, z=\infty \right) .  \end{eqnarray*} 
To be able to use the residue theorem, this time we have to consider the function
  \[ - \frac{ 1 }{ z^{ a_{ 1 } } - 1 } \cdots \frac{ 1 }{ z^{ a_{ n } } - 1 } \ \frac{ 1 }{ z - 1 } \ z^{ t - 1 } = (-1)^{ n } F_{ t } (z) \]
The residues at the finite poles of $F_{ t } $ can be computed as before, with $t$ replaced by $-t$, and the proof of 
the following theorem is completely analogous to Theorem \ref{closure}: 
\begin{theorem} Let ${\cal P}$ be given by {\rm (\ref{polytope})}, with $a_{ 1 } , \dots, a_{ n } $ pairwise relatively prime. Then 
  \[ L ( {\cal P}^{\circ} , t ) = (-1)^{ n } R_{ t } \left( a_{ 1 } , \dots, a_{ n } \right) + \sum_{ j=1 }^{ n } \sigma_{t} ( a_{ 1 } , \dots, {\hat a_{ j } } , \dots , a_{ n } , 1 ; a_{ j } )  \] 
\end{theorem} 
As an immediate consequence we get the remarkable
\begin{corollary}[Ehrhart-Macdonald Reciprocity Law]\label{ehr} 
  \[ L ( {\cal P}^{\circ} , -t ) = (-1)^{ n } L ( \overline{\cal P} , t ) \ . \]
\end{corollary}
This result was conjectured for convex rational polytopes by Ehrhart \cite{ehrhart}, 
and first proved by Macdonald \cite{macdonald}.

Of particular interest is the number of lattice points on the boundary of $t {\cal P}$. 
Besides computing $ L ({\cal P}^{\circ} , t )$ and $ L ( \overline{\cal P} , t )$ and taking 
differences, we can also adjust our method to this situation, especially if we are 
interested in only \emph{parts} of the boundary. 
As an example, we will compute $p_A' (t) $ 
as defined in the introduction (\ref{N_t}), which appears in the context of the Frobenius problem.  
Again, for reasons of simplicity we assume in the following that $a_{1}, \dots a_{n}$ are 
\emph{pairwise coprime} positive integers. 

This time we interpret 
  \[ p_A' (t) = \# \left\{ (m_{1}, \dots, m_{n}) \in \N^{n} : \sum_{k=1}^{n} m_{k} a_{k} = t \right\} \]
as the Taylor coefficient of $z^{ t }$ of the function
\begin{eqnarray*} &\mbox{}& \left( z^{ a_{ 1 } } + z^{ 2 a_{ 1 } } + \dots  \right) \cdots \left( z^{ a_{ n } } + z^{ 2 a_{ n } } + \dots  \right) \\
                  &\mbox{}& \quad = \frac{ z^{ a_{ 1 } } }{ 1 - z^{ a_{ 1 } } } \cdots \frac{ z^{ a_{ n } } }{ 1 - z^{ a_{ n } } } \ . \end{eqnarray*} 
That is, 
\begin{eqnarray*} &\mbox{}& p_A' (t) = \mbox{Res} \left( \frac{ z^{ a_{ 1 } } }{ 1 - z^{ a_{ 1 } } } \cdots \frac{ z^{ a_{ n } } }{ 1 - z^{ a_{ n } } } \ z^{ - t - 1 } , z=0 \right)  \\
                  &\mbox{}& \quad = \mbox{Res} \left( \frac{ -1 }{ z^{ 2 }  } \ \frac{ 1 }{ z^{ a_{ 1 } } - 1 } \cdots \frac{ 1 }{ z^{ a_{ n } } - 1 } \ z^{ t + 1 }, z=\infty \right) .  \end{eqnarray*} 
Thus, we have to find the other residues of 
  \[ G_{ t } (z) := \frac{ z^{ t - 1 } }{ \left( z^{ a_{ 1 } } - 1 \right) \cdots \left( z^{ a_{ n } } - 1 \right) } = \left( z-1 \right) F_{ t }(z)  \ , \]
since 
  \begin{equation}\label{N_tRes} p_A' (t) = - \mbox{ Res} \left( G_{ t } (z) , z = \infty \right) . \end{equation}
$G_{ t } $ has its other poles at all $a_{1}^{\text{th}}, \dots, a_{n}^{\text{th}}$ roots of unity. 
Again, note that $G_{t}$ has \emph{simple} poles at all the nontrivial roots of unity. Let 
$\lambda$ be a nontrivial $a_{1}^{\text{th}}$ root of unity, then
\begin{eqnarray*} &\mbox{}& \mbox{Res} \left( G_{ t } (z) , z = \lambda \right) = \frac{ \lambda^{ t - 1 } }{ \left( \lambda^{ a_{ 2 } } - 1 \right) \cdots \left( \lambda^{ a_{ n } } - 1 \right) } \ \mbox{Res} \left( \frac{ 1 }{ z^{ a_{ 1 } } - 1 } , z = \lambda \right) \\
                  &\mbox{}& \quad = \frac{ \lambda^{ t } }{ a_{1} \left( \lambda^{ a_{ 2 } } - 1 \right) \cdots \left( \lambda^{ a_{ n } } - 1 \right) } . \end{eqnarray*} 
Adding up all the nontrivial $a_{1}^{\text{th}}$ roots of unity, we obtain
\begin{eqnarray*} &\mbox{}& \sum_{ \lambda^{ a_{1} } = 1 \not= \lambda } \mbox{Res} \left( G_{ t } (z) , z = \lambda \right) = \frac{1}{a_{1}} \sum_{ \lambda^{ a_{1} } = 1 \not= \lambda } \frac{ \lambda^{ t } }{ \left( \lambda^{ a_{ 2 } } - 1 \right) \cdots \left( \lambda^{ a_{ n } } - 1 \right) } \\
                  &\mbox{}& \quad = \sigma_{t} \left(a_{ 2 } , \dots, a_{ n } ; a_{ 1 }  \right) . \end{eqnarray*} 
Together with the similar residues at the other roots of unity, (\ref{N_tRes}) gives us
\begin{theorem}\label{frobnr}
  \[ p_A' (t) = R'_{ t }\left( a_{ 1 } , \dots, a_{ n } \right) + \sum_{ j=1 }^{ n } \sigma_{t} ( a_{ 1 } , \dots, {\hat a_{ j } } , \dots , a_{ n } ; a_{ j } ) \ , \] 
where $ R'_{ t } \left( a_{ 1 } , \dots, a_{ n } \right) = \mbox{\rm Res} \left( G_{ t } (z) , z=1 \right) $. 
\end{theorem}

$R'$ is as easily computed as before, the first values are
\begin{eqnarray*} &\mbox{}& R'_{ t }(a,b) = \frac{ t }{ ab } - \frac{ 1 }{ 2 } \left( \frac{ 1 }{ a } + \frac{ 1 }{ b }  \right) \\
                  &\mbox{}& R'_{ t }(a,b,c) = \frac{ t^{ 2 }  }{ 2abc } - \frac{ t }{ 2 } \left( \frac{ 1 }{ ab } + \frac{ 1 }{ ac } + \frac{ 1 }{ bc }  \right) + \\
                  &\mbox{}& \quad + \frac{ 1 }{ 12 } \left( \frac{ 3 }{ a } + \frac{ 3 }{ b } + \frac{ 3 }{ c } + \frac{ a }{ bc } + \frac{ b }{ ac } + \frac{ c }{ ab }  \right) \\
                  &\mbox{}& R'_{ t }(a,b,c,d) = \frac{ t^{ 3 }  }{ 6abcd } - \frac{ t^{ 2 }  }{ 4 } \left( \frac{ 1 }{ abc } + \frac{ 1 }{ abd } + \frac{ 1 }{ acd } + \frac{ 1 }{ bcd }  \right) \\
                  &\mbox{}& \qquad + \frac{ t }{ 12 } \left( \frac{ 3 }{ ab } + \frac{ 3 }{ ac } + \frac{ 3 }{ ad } + \frac{ 3 }{ bc } + \frac{ 3 }{ bd } + \frac{ 3 }{ cd } + \frac{ a }{ bcd } + \frac{ b }{ acd } + \frac{ c }{ abd } + \frac{ d }{ abc }  \right) \\
                  &\mbox{}& \qquad - \frac{ 1 }{ 24 } \left( \frac{ a }{ bc } + \frac{ a }{ bd } + \frac{ a }{ cd } + \frac{ b }{ ad } + \frac{ b }{ ac } + \frac{ b }{ cd } + \frac{ c }{ ab } + \frac{ c }{ ad } + \frac{ c }{ bd } + \frac{ d }{ ab } + \frac{ d }{ ac } + \frac{ d }{ bc } \right) \\ 
                  &\mbox{}& \qquad - \frac{ 1 }{ 8 } \left( \frac{ 1 }{ a } + \frac{ 1 }{ b } + \frac{ 1 }{ c } + \frac{ 1 }{ d } \right) \ . \end{eqnarray*} 
A general formula for $ R'_{ t } \left( a_{ 1 } , \dots, a_{ n } \right) $ was recently discovered in \cite{bgk}. 

For generalizations, note that we can apply our method to any tetrahedron given in the form 
(\ref{polytope}), with the $ a_{ k } $'s replaced by any rational numbers. Moreover, any convex 
rational polytope (that is, a convex polytope whose vertices have rational coordinates) 
can be described by a finite number of inequalities over the rationals. In other words, a 
convex lattice polytope ${\cal P} $ is an intersection of finitely many half-spaces. 
This description of the polytope leads to an integral in several complex variables, as discussed 
in \cite[Theorem 8]{beck} for lattice polytopes. 


\section{\normalsize The Fourier method}\label{sinaiSec}
In this section we outline a Fourier-analytic method that achieves the same results. 
Although the theory is a little harder, the method is of independent interest. 
It draws connections to Brion's theorem on generating functions \cite{brion} and to 
the basic results of \cite{diaz}. 

To be concrete, we illustrate the general case with the 2-dimensional rational triangle 
${\cal P} $ whose vertices are $ v_{ 0 } = (0,0) $, $ v_{ 1 } = \left( \frac{ t }{ a } , 0 \right)  $, and 
$ v_{ 2 } = \left( 0 , \frac{ t }{ b }  \right)  $. As before, the number of lattice points in the 1-dimensional 
hypotenuse of this right triangle is 
  \[ p_{ \{ a, b \} }' (t) = \# \left\{ (m,n) \in \N^{ 2 } : am + bn = t \right\} . \] 
We denote the tangent cone to ${\cal P} $ at the vertex $ v_{ i } $ by $ K_{ i } $. 
We recall that the exponential sum attached to the cone $K$ (with vertex $v$) is by definition 
  \begin{equation}\label{cone} \sigma_{ K } (s) = \sum_{ m \in \Z^{ n } \cap K } e^{ -2 \pi \left< s,m \right>  } \ , \end{equation} 
where $s$ is any complex vector that makes the infinite sum (\ref{cone}) converge. 
An equivalent formulation of (\ref{cone}) which appears more combinatorial is
  \begin{equation}\label{cone2} \sigma_{ K } (x) = \sum_{ m \in \Z^{ n } \cap K } x^{ m }  \ , \end{equation} 
where $ x^{ m } = x_{ 1 }^{ m_{ 1 }  }  \cdots x_{ n }^{ m_{ n }  } $ and $ x_{ j } = e^{ - 2 \pi s_{ j }  } $. 

In general dimension, let the vertices of the rational polytope ${\cal P} $ be $ v_{ 1 } , \dots , v_{ l } $. 
Let the corresponding tangent cone at $ v_{ j } $ be $ K_{ j } $. Finally, let the \emph{finite} 
exponential sum over ${\cal P} $ be 
  \begin{equation}\label{cone3} \sigma_{\cal P} (s) = \sum_{ m \in \Z^{ n } \cap {\cal P} } e^{ -2 \pi \left< s,m \right>  } \ . \end{equation} 
Then there is the basic result that each exponential sum (\ref{cone2}) is a rational function 
of $x$, and the following theorem relates these rational functions \cite{brion}: 
\begin{theorem}[Brion]\label{brionth} 
For a generic value of $ s \in \C^{ n } $, 
  \begin{equation}\label{brioneq} \sigma_{\cal P} (s) = \sum_{ i=1 }^{ l } \sigma_{ K_{ i }  } (s) \ . \end{equation} 
\end{theorem}
This result allows us to transfer the enumeration of lattice points in ${\cal P} $ to the 
enumeration of lattice points in the tangent cones $ K_{ i } $ at the vertices of ${\cal P} $, 
an easier task. In the theorem above, `generic value of $s$' means any $ s \in \C^{ n } $ 
for which these rational functions do not blow up to infinity. 

To apply these results to our given rational triangle ${\cal P} $, 
we first employ the methods of \cite{diaz} to get an explicit formula for the exponential 
sum for each tangent cone of ${\cal P} $. Then, by Brion's theorem on 
tangent cones, the sum of the three exponential sums attached to the tangent cones equals the 
exponential sum over ${\cal P} $. Canceling the singularities arising from each tangent cone, 
and letting $ s \to 1 $, we get the explicit formula of the previous section for the number 
of lattice points in the rational triangle ${\cal P} $. 

In our case, $ K_{ 1 } $ is generated by the two rational vectors $ -v_{ 1 } $ and $ v_{ 2 } - v_{ 1 } $. 
We form the matrix 
  \[ A_{ 1 } = \left( \begin{array}{cc} - \frac{ t }{ a }  & - \frac{ t }{ a } \\ 
                                                   0       & \frac{ t }{ b }   \end{array} \right) \ , \]
whose columns are the vectors that generate the cone $ K_{ 1 } $. Once we compute $ \sigma_{ K_{ 1 }  } (s) $, 
$ \sigma_{ K_{ 2 }  } (s) $ will follow by symmetry. The easiest exponential sum to compute is 
\begin{eqnarray*} &\mbox{}& \sigma_{ K_{ 0 }  } (s) = \sum_{ m \in \Z^{ 2 } \cap K_{ 0 }  } e^{ -2 \pi \left< s,m \right>  } 
                                                    = \sum_{ {m_{ 1 } \geq 0} \atop {m_{ 2 } \geq 0} } e^{ -2 \pi ( m_{ 1 } s_{ 1 } + m_{ 2 } s_{ 2 }  ) } \\ 
                  &\mbox{}& \qquad \qquad = \frac{ 1 }{ \left( 1 - e^{ -2 \pi s_{ 1 }  } \right) \left( 1 - e^{ -2 \pi s_{ 2 }  } \right) } \ . \end{eqnarray*} 
To compute $ \sigma_{ K_{ i }  } (s) $ ($ i \not= 0 $), we first translate the cone $K_{ i }$ 
by the vector $ - v_{ i }  $ so that its new vertex is the origin. We therefore let $ K = K_{ i } - v_{ i } $, 
and the following elementary lemma illustrates how a translation affects the Fourier transform. 
Let 
  \[ \chi_K (x) = \left\{ \begin{array}{ll} 1 & \text{ if } x \in K , \\ 
                                            0 & \text{ if } x \not\in K \end{array} \right. \] 
denote the characteristic function of $K$. 
\begin{lemma} Let 
  \[ F_{ v } (x) = \chi_{ K + v } (x) \ e^{ -2 \pi \left< s,m \right>  } \] 
for $ x \in \R^{ n } , s \in \C^{ n } $. Then 
  \[ { \hat F }_{ v } ( \xi ) = { \hat \chi }_{ K } ( \xi + i s ) \ e^{ -2 \pi i \left< \xi + i s , v \right>  } \] 
\end{lemma} 

\begin{proof} 
\begin{eqnarray*} &\mbox{}& { \hat F }_{ v } ( \xi ) = \int_{ \R^{ n }  } \chi_{ K + v } (x) \ e^{ -2 \pi \left< s,m \right>  } \ e^{ 2 \pi i \left< \xi , x \right>  } \ dx \\ 
                  &\mbox{}& \quad = \int_{ \R^{ n }  } e^{ 2 \pi i \left< \xi + i s , x \right>  } \ \chi_{ K + v } (x) \ dx \\ 
                  &\mbox{}& \quad = \int_{ \R^{ n }  } e^{ 2 \pi i \left< \xi + i s , y-v \right>  } \ \chi_{ K } (y) \ dy \\ 
                  &\mbox{}& \quad =  e^{ -2 \pi i \left< \xi + i s , v \right>  } \int_{ \R^{ n }  } \ e^{ 2 \pi i \left< \xi + i s , y \right>  } \ \chi_{ K } (y) \ dy \\ 
                  &\mbox{}& \quad = e^{ -2 \pi i \left< \xi + i s , v \right>  } \ { \hat \chi }_{ K } ( \xi + i s ) \end{eqnarray*} 
\end{proof} 

This lemma also shows why it is useful to study the Fourier transform of $K$ at \emph{complex} 
values of the variable; that is, at $ \xi + i s $. We study $ F(x) $ because (\ref{cone}) can 
be rewritten as 
  \[ \sigma_{ K_{ 0 } + v } (s) = \sum_{ m \in \Z^{ n } } \chi_{ K_{ 0 } + v } \ e^{ -2 \pi \left< s,m \right>  } = \sum_{ m \in \Z^{ n } } F_{ v } (m) \ . \]
All of the lemmas of \cite{diaz} remain true in this rational polytope context. The idea is to 
apply Poisson summation to $ \sum_{ m \in \Z^{ n } } F_{ v } (m) $ and write formally 
  \[ \sum_{ m \in \Z^{ n } } F_{ v } (m) = \sum_{ m \in \Z^{ n } } { \hat F }_{ v } (m)  \]
The right-hand side diverges, though, and some smoothing completes the picture. Because the 
steps are identical to those in \cite{diaz}, we omit the ensuing details. Let $ \xi_{ a } = e^{ \frac{ 2 \pi i }{ a }  } $. 
We get 
\begin{eqnarray}\label{coth} &\mbox{}& \sigma_{ K_{ 1 }  } \left( s_{ 1 } , s_{ 2 }  \right) = \frac{ \xi_{ a } ^{ t s_{ 1 }  }  }{ 4 a } \sum_{ r=0 }^{ a-1 } \xi_{ a }^{ rt } \left( \coth \frac{ \pi b }{ t } \left( s_{ 1,2 } + \frac{ irt }{ a }  \right) - 1 \right) \\ 
                             &\mbox{}& \qquad \qquad \qquad \qquad \qquad \cdot \left( \coth \frac{ \pi }{ t } \left( s_{ 1,1 } + \frac{ irt }{ a }  \right) + 1 \right) \ , \nonumber \end{eqnarray} 
where 
 \[ s_{ 1,1 } = \left< s, \mbox{ generator 1 of } K_{ 1 }  \right> = \left< \left( s_{ 1 } , s_{ 2 }  \right) , \left( - \frac{ t }{ a } , 0 \right)  \right> = - \frac{ t s_{ 1 }  }{ a } \] 
and 
  \[ s_{ 1,2 } = \left< s, \mbox{ generator 2 of } K_{ 1 }  \right> = \left< \left( s_{ 1 } , s_{ 2 }  \right) , \left( - \frac{ t }{ a } , \frac{ t }{ b }  \right)  \right> = - \frac{ t s_{ 1 }  }{ a } + \frac{ t s_{ 2 }  }{ b } \ . \] 
By (\ref{brioneq}), we have 
  \[ \# \left\{ \Z^{ 2 } \cap t {\cal P} \right\} \ = \ \sum_{ m \in \Z^{ 2 } \cap t {\cal P}  } 1 \ = \ \lim_{ s \to 0 } \left( \sigma_{ K_{ 0 }  } (s) + \sigma_{ K_{ 1 }  } (s) + \sigma_{ K_{ 2 }  } (s) \right) \ . \]
Using the explicit description of $ \sigma_{ K_{ i }  } (s) $ in terms of cotangent 
functions, we can cancel their singularities at $ s=0 $ and simply add the holomorphic 
contributions to $ \sigma_{ K_{ i }  } (s) $ at $ s=0 $. The left-hand side of (\ref{brioneq}) 
is holomorphic in $s$, so that we are guaranteed that the singularities on the right-hand side 
cancel each other. 

The only term in the finite sum (\ref{coth}) that contributes a singularity at $s=0$ is the 
$ r=0 $ term. We expand the three exponential sums $ \sigma_{ K_{ i }  } (s) $ into their 
Laurent expansions about $ s=0 $. Here we only require the first 3 terms of their Laurent 
expansions. In dimension $n$ we would require the first $ n+1 $ terms; otherwise every step 
is the same in general dimension $n$. 

We make use of the Laurent series
  \[ \frac{ 1 }{ 1 - e^{ - \alpha s }  } = \frac{ 1 }{ \alpha s } + \frac{ 1 }{ 2 } + \frac{ \alpha s }{ 12 } + O \left( s^{ 2 }  \right)  \]
near $ s=0 $, as well as the Laurent series for $ \cot \pi s $ near $ s=0 $. After expanding 
each cotangent in (\ref{coth}) for $ \sigma_{ K_{ 0 }  } (s) $, $ \sigma_{ K_{ 1 }  } (s) $ 
and $ \sigma_{ K_{ 2 }  } (s) $ and letting $ s \to 0 $, we obtain Theorem \ref{closure} above as 
\begin{eqnarray*} &\mbox{}& L ( \overline{\cal P} , t ) = \frac{ t^{ 2 }  }{ 2ab } + \frac{ t }{ 2 } \left( \frac{ 1 }{ a } + \frac{ 1 }{ b } + \frac{ 1 }{ ab }  \right) \\ 
                  &\mbox{}& \qquad + \frac{ 1 }{ 4 } \left( 1 + \frac{ 1 }{ a } + \frac{ 1 }{ b }  \right) + \frac{ 1 }{ 12 } \left( \frac{ a }{ b } + \frac{ b }{ a } + \frac{ 1 }{ ab }  \right) \\ 
                  &\mbox{}& \qquad + \frac{ 1 }{ a } \sum_{ r=1 }^{ a-1 } \frac{ \xi_{ a }^{ rt }   }{ \left( 1 - \xi_{ a }^{ rb }   \right) \left( 1 - \xi_{ a }^{ r }   \right) } + \frac{ 1 }{ b } \sum_{ r=1 }^{ b-1 } \frac{ \xi_{ b }^{ rt }   }{ \left( 1 - \xi_{ b }^{ ra }   \right) \left( 1 - \xi_{ b }^{ r }   \right) } \ . \end{eqnarray*} 
Note that, as before, the periodic portion of $ L ( \overline{\cal P} , t ) $ is entirely 
contained in the ``constant'' $t$ term. 
By Ehrhart's reciprocity law (Corollary \ref{ehr}, \cite{ehrhart}), there is a similar 
expression for $ L ( {\cal P} , t ) $, and taking 
  \[ L ( \overline{\cal P} , t ) - L ( {\cal P} , t ) - \left[ \frac{ t }{ a }  \right] - \left[ \frac{ t }{ b }  \right] - 1 \] 
gives us $ p_{ \{ a, b \} } (t) $. The same analysis gives us Theorem \ref{closure} in $ \R^{ n } $. 


\section{\normalsize The Fourier-Dedekind sum}\label{FDSum}
In the derivation of the various lattice count formulas, we naturally arrived at 
the Fourier-Dedekind sum
  \[ \sigma_{t} \left( c_{ 1 } , \dots, c_{ n } ; c  \right) = \frac{1}{c} \sum_{ \lambda^{ c } = 1 \not= \lambda } \frac{ \lambda^{ t } }{ \left( \lambda^{ c_{ 1 } } - 1 \right) \cdots \left( \lambda^{ c_{ n } } - 1 \right) } \ . \]
This expression is a generalization of the classical Dedekind sum $\s(h,k)$ \cite{grosswald} 
and its various generalizations \cite{dieter,gessel,meyer,rademacher,zagier}. 
In fact, an easy calculation shows
\begin{eqnarray*} &\mbox{}& \sigma_{ 0 }  \left( a, 1; c \right) = \frac{ 1 }{ c } \sum_{ \lambda^{ c } = 1 \not= \lambda } \frac{ 1 }{ \left( \lambda^{ a } - 1 \right) \left( \lambda - 1 \right) } = \\
                  &\mbox{}& \quad = \frac{ 1 }{ 4 } - \frac{ 1 }{ 4c } - \frac{ 1 }{ 4c } \sum_{ k=1 }^{ c-1 } \cot \frac{ \pi k a }{ c } \cot \frac{ \pi k }{ c } = \frac{ 1 }{ 4 } - \frac{ 1 }{ 4c } - \s (a,c) \ . \end{eqnarray*}
In general, note that $\sigma_{t} \left(c_{ 1 } , \dots, c_{ n } ; c  \right)$ is a rational 
number: It is an element of the cyclotomic field of $c^{\text{th}}$ roots of unity, and invariant under all 
Galois transformations of this field. 

Some obvious properties are
\begin{eqnarray} &\mbox{}& \sigma_{t} \left( c_{ 1 } , \dots, c_{ n } ; c  \right) = \sigma_{t} \left( c_{ \pi (1) } , \dots, c_{ \pi (n) } ; c  \right) \quad \mbox{ for any } \pi \in S_{ n } \nonumber \\
                 &\mbox{}& \sigma_{t} \left( c_{ 1 } , \dots, c_{ n } ; c  \right) = \sigma_{(t \mbox{ \scriptsize mod } c) } \left( c_{ 1 } \mbox{ mod } c , \dots, c_{ n } \mbox{ mod } c ; c  \right) \label{period} \\
                 &\mbox{}& \sigma_{t} \left( c_{ 1 } , \dots, c_{ n } ; c  \right) = \sigma_{bt} \left( b c_{ 1 } , \dots, b c_{ n } ; c  \right) \quad \mbox{ for any } b \in \Z \mbox{ with } (b,c)=1 \nonumber 
\end{eqnarray} 
We can get more familiar-looking formulas for $\sigma_{t}$ in certain dimensions. For example, counting points in dimension 1, we find that
  \[ L ( \overline{\cal P} , t ) = \# \left\{ m \in \Z : m \geq 0 , m c \leq t \right\} = \left\lfloor \frac{ t }{ c } \right\rfloor + 1 \ , \]
so that Theorem \ref{closure} implies 
  \begin{equation}\label{fouriersingle} \sigma_{-t} ( 1 ; c ) = \frac{1}{c} \sum_{ \lambda^{ c } = 1 \not= \lambda } \frac{ \lambda^{ -t } }{ \left( \lambda - 1 \right) } = \frac{ t }{ c } - \left\lfloor \frac{ t }{ c } \right\rfloor - \frac{ 1 }{ 2 }  + \frac{ 1 }{ 2c } = \left( \left( \frac{t}{c} \right) \right) + \frac{1}{2c} \ .  \end{equation}
Here, $ (( x )) = x - \lfloor x \rfloor - 1/2 $ is a sawtooth function (differing slightly from the one appearing in the classical 
Dedekind sums). This restates the well-known finite Fourier expansion of the sawtooth function (see, e.g., \cite{grosswald}). 

As another example, we reformulate 
  \[ \sigma_{t} ( a, b; c ) = \frac{1}{c} \sum_{ \lambda^{ c } = 1 \not= \lambda } \frac{ \lambda^{ t } }{ \left( \lambda^{ a } - 1 \right) \left( \lambda^{ b } - 1 \right) } \]
by means of finite Fourier series. Consider
\begin{eqnarray} &\mbox{}& \sigma_{t} ( a; c ) = \frac{1}{c} \sum_{ \lambda^{ c } = 1 \not= \lambda } \frac{ \lambda^{ -t } }{ \left( \lambda^{ a }  - 1 \right) } = \frac{1}{c} \sum_{ k=1 }^{ c - 1 } \frac{ \xi^{ k t } }{ \left( \xi^{ k a } - 1 \right) } = \frac{1}{c} \sum_{ k=1 }^{ c - 1 } \frac{ \xi^{ k a^{-1} t } }{ \left( \xi^{ k } - 1 \right) } \nonumber \\ 
                 &\mbox{}& \quad = \left( \left( \frac{ -a^{-1} t }{ c } \right) \right) + \frac{1}{2c} \ , \label{fourierhelp} \end{eqnarray} 
where $\xi$ is a primitive $c^{\text{th}}$ root of unity and $ a a^{-1} \equiv 1 $ mod $c$; 
here, the last equality follows from (\ref{fouriersingle}). 
We use the well-known convolution theorem for finite Fourier series: 
\begin{theorem} Let $ f(t) = \frac{ 1 }{ N }  \sum_{ k=0 }^{ N-1 } a_{ k } \xi^{ kt } $ and 
$ g(t) = \frac{ 1 }{ N } \sum_{ k=0 }^{ N-1 } b_{ k } \xi^{ kt } $, 
where $\xi$ is a primitive $N^{\text{th}}$ root of unity. Then 
  \[ \frac{ 1 }{ N } \sum_{ k=0 }^{ N-1 } a_{ k } b_{ k } \xi^{ kt } = \sum_{ m=0 }^{ N-1 } f(t-m) g(m) \ .  \] 
\end{theorem} 
Hence by (\ref{fourierhelp}), 
\begin{eqnarray*} &\mbox{}& \sigma_{t} ( a, b; c ) = \sum_{ m=0 }^{ c-1 } \sigma_{t-m} ( a; c ) \sigma_{m} ( b; c ) \\
                  &\mbox{}& \qquad = \sum_{ m=0 }^{ c-1 } \left[ \left( \left( \frac{ -a^{-1} (t-m) }{ c } \right) \right) + \frac{1}{2c} \right] \left[ \left( \left( \frac{ -b^{-1} m }{ c } \right) \right) + \frac{1}{2c} \right] \\
                  &\mbox{}& \qquad = \sum_{ m=0 }^{ c-1 } \left( \left( \frac{ a^{-1} (m-t) }{ c } \right) \right) \left( \left( \frac{ -b^{-1} m }{ c } \right) \right) - \frac{1}{4c} \ . \end{eqnarray*}
Here, $ a a^{-1} \equiv b b^{-1} \equiv 1 $ mod $c$. The last equality follows from 
\[ \sum_{ m=0 }^{ c-1 } \left( \left( \frac{ m }{ c } \right) \right) = - \frac{ 1 }{ 2 } \ . \] 
Furthermore, by the periodicity of $ ((x)) $, 
\begin{equation}\label{fouriertriple} \sigma_{t} ( a, b; c ) = \sum_{ m=0 }^{ c-1 } \left( \left( \frac{ - a^{-1} ( b m + t ) }{ c } \right) \right) \left( \left( \frac{ m }{ c } \right) \right) - \frac{1}{4c} \ . \end{equation}
The expression on the right is, up to a trivial term, a special case of a 
\emph{Dedekind-Rademacher sum} \cite{dieter,knuth,meyer,rademacher}. 
It is a curious fact that the function $ \sigma_{t} ( a, b; c ) $ is the nontrivial part of a 
multiplier system of a weight-0 modular form \cite[p.~121]{robins}. 

We conlude this section by proving two reciprocity laws for Fourier-Dedekind sums. The first one is 
equivalent to Zagier's reciprocity law for his \emph{higher dimensional Dedekind sums} \cite{zagier}. 
They are essentially Fourier-Dedekind sums with $t=0$, that is, trivial numerators. 
\begin{theorem}\label{zag} For pairwise relatively prime integers $ a_{1} , \dots , a_{n} $, 
  \[ \sum_{ j=1 }^{ n } \sigma_{0} ( a_{ 1 } , \dots, {\hat a_{ j } } , \dots , a_{ n } ; a_{ j } ) = 1 - R'_{ 0 }\left( a_{ 1 } , \dots, a_{ n } \right) \ , \]
where $ R'_{ t } $ is the rational function given in Theorem \ref{frobnr}. 
\end{theorem} 
\begin{proof} 
It is well known \cite{ehrhart} that the constant term of a \emph{lattice} polytope (that is, a polytope with integral vertices) equals the Euler characteristic of the polytope. 
Consider the polytope 
\[ \left\{ (x_{ 1 } , \dots , x_{ n }) \in \R_{ >0 }^{ n }  : \sum_{k=1}^{n} x_{k} a_{k} = 1 \right\} , \] 
whose dilates correspond to the quantor $p_A' (t) $ of Theorem \ref{frobnr}. 
If we dilate this polytope only by multiples of $ a_{1} \cdots a_{n} $, say $ t = a_{1} \cdots a_{n} w $, 
we obtain the dilates of a lattice polytope. Theorem \ref{frobnr} simplifies for these $t$ to 
  \[ p_A' ( a_{1} \cdots a_{n} w ) = R'_{ a_{1} \cdots a_{n} w } \left( a_{ 1 } , \dots, a_{ n } \right) + \sum_{ j=1 }^{ n } \sigma_{0} ( a_{ 1 } , \dots, {\hat a_{ j } } , \dots , a_{ n } ; a_{ j } ) \ , \] 
using the periodicity of $ \sigma_{t} $ (\ref{period}). 
On the other hand, we know that the constant term (in terms of $w$) is the Euler characteristic of 
the polytope and hence equals 1, which yields the identity 
  \[ 1 = R'_{ 0 } \left( a_{ 1 } , \dots, a_{ n } \right) + \sum_{ j=1 }^{ n } \sigma_{0} ( a_{ 1 } , \dots, {\hat a_{ j } } , \dots , a_{ n } ; a_{ j } ) \ . \] 
\end{proof} 

The second one is a new reciprocity law, which generalizes the following \cite{gessel} 
\begin{theorem}[Gessel]\label{ges} Let $m$ and $n$ be relatively prime and suppose that $ 0 \leq r < m+n $. Then 
  \begin{eqnarray*} &\mbox{}& \frac{ 1 }{ m }  \sum_{ \lambda^{ m } = 1 \not= \lambda } \frac{ \lambda^{ r+1 }  }{ \left( \lambda^{ n } - 1 \right) \left( \lambda - 1 \right)  } + \frac{ 1 }{ n }  \sum_{ \lambda^{ n } = 1 \not= \lambda } \frac{ \lambda^{ r+1 }  }{ \left( \lambda^{ m } - 1 \right) \left( \lambda - 1 \right)  } \\ 
                    &\mbox{}& \qquad = - \frac{1}{12} \left( \frac{m}{n} + \frac{n}{m} + \frac{1}{mn} \right) + \frac{1}{4} \left( \frac{1}{m} + \frac{1}{n} - 1 \right) \\ 
                    &\mbox{}& \qquad \qquad + \frac{r}{2} \left( \frac{1}{m} + \frac{1}{n} - \frac{1}{mn} \right) - \frac{ r^{2} }{ 2mn } \ . \end{eqnarray*} 
\end{theorem} 
It is not hard to see that Gessel's theorem follows as the two-dimensional case of 
\begin{theorem} Let $ a_{1} , \dots , a_{n} $ be pairwise relatively prime integers and $0 < t < a_{ 1 } + \dots + a_{ n }$. Then 
  \[ \sum_{ j=1 }^{ n } \sigma_{t} ( a_{ 1 } , \dots, {\hat a_{ j } } , \dots , a_{ n } ; a_{ j } )  = - R'_{ t } \left( a_{ 1 } , \dots, a_{ n } \right)  \ , \]
where $ R'_{ t } $ is the rational function given in Theorem \ref{frobnr}. 
\end{theorem} 
\begin{proof} 
By definition, $p_A' (t) = 0$ if $0 < t < a_{ 1 } + \dots + a_{ n }$. 
Hence Theorem \ref{frobnr} yields an identiy for these values of $t$: 
  \[ 0 = R'_{ t } \left( a_{ 1 } , \dots, a_{ n } \right) + \sum_{ j=1 }^{ n } \sigma_{t} ( a_{ 1 } , \dots, {\hat a_{ j } } , \dots , a_{ n } ; a_{ j } ) \ . \] 
\end{proof} 

It is worth noticing that both Theorems \ref{zag} and \ref{ges} imply the reciprocity law for the classical 
Dedekind sum $ \s (a,b) $. 
It should be finally mentioned that in special cases there are other reciprocity laws, 
for example, for the sum appearing on the right-hand side in (\ref{fouriertriple}) \cite{dieter,rademacher}. 
We note that, as a consequence, we can compute $ \sigma_{t} ( a, b; c ) $ in polynomial time. 


\section{\normalsize The Frobenius problem}\label{lastsec}
In this last section we apply Theorem \ref{frobnr} (the explicit formula for 
$ p_A' (t) $) to Frobenius's original problem. 
As an example, we will discuss the 3-dimensional case. 
Note that a bound for dimension 3 yields a bound for the general case: It can be 
easily verified that
\begin{equation}\label{threegen} f(a_{1}, \dots, a_{n}) \leq f(a_{1}, a_{ 2 } , a_{3}) + a_{ 4 } + \dots + a_{ n }  \end{equation}
Furthermore, in dimension 3 it suffices to assume that $ a_{ 1 } , a_{ 2 } , a_{ 3 } $ are pairwise coprime, 
due to Johnson's formula \cite{johnson}: If $ g = ( a_{ 1 } , a_{ 2 } ) $, then 
\begin{equation}\label{coprimegen} f(a_{1}, a_{ 2 } , a_{3}) = g \cdot f \left( \frac{ a_{ 1 }  }{ g } , \frac{ a_{ 2 }  }{ g } , a_{ 3 }  \right) \ . \end{equation}


Now assume $ a, b, c $ pairwise relatively prime, and recall (\ref{fouriertriple}): 
\[ \sigma_{t} ( a, b; c ) = \sum_{ m=0 }^{ c-1 } \left( \left( \frac{ - a^{-1} ( b m + t ) }{ c } \right) \right) \left( \left( \frac{ m }{ c } \right) \right) - \frac{1}{4c} \ , \]
where $ a a^{-1} \equiv 1 $ mod $c$. 
We will use the Cauchy-Schwartz inequality 
  \begin{equation}\label{perm} \left| \sum_{ k=1 }^{ n } a_{ k } a_{ \pi(k) }  \right| \leq \sum_{ k=1 }^{ n } a_{ k }^{ 2 }  \ . \end{equation} 
Here $ a_{ k } \in \R $, and $\pi \in S_{ n } $ is a permutation. 
Since $ ( a^{-1} b , c ) = 1 $, we can use (\ref{perm}) to obtain
\begin{eqnarray*} &\mbox{}& \sigma_{t} ( a, b; c ) \geq - \sum_{ m=0 }^{ c-1 } \left( \left( \frac{ m }{ c } \right) \right)^{ 2 }  - \frac{1}{4c} = \sum_{ m=0 }^{ c-1 } \left( \frac{ m }{ c } - \frac{ 1 }{ 2 } \right)^{ 2 }  - \frac{1}{4c} \\
                  &\mbox{}& \qquad = - \frac{ 1 }{ c^{ 2 }  } \frac{ (2c-1) (c-1) c }{ 6 } + \frac{ 1 }{ c } \frac{ c (c-1) }{ 2 } - \frac{ c }{ 4 } - \frac{1}{4c} \\
                  &\mbox{}& \qquad = - \frac{ c }{ 12 } - \frac{ 1 }{ 12 c }  \ . \end{eqnarray*}
This also restates Rademacher's bound on the classical Dedekind sums \cite{grosswald}. 
Using this in the formula for dimension 3 (remark after Theorem \ref{frobnr}), we get
\begin{eqnarray*} &\mbox{}& p_{\{a,b,c\}}' (t) \geq \frac{ t^{ 2 }  }{ 2abc } - \frac{ t }{ 2 } \left( \frac{ 1 }{ ab } + \frac{ 1 }{ ac } + \frac{ 1 }{ bc }  \right) \\
                  &\mbox{}& \qquad \qquad \qquad + \frac{ 1 }{ 12 } \left( \frac{ 3 }{ a } + \frac{ 3 }{ b } + \frac{ 3 }{ c } + \frac{ a }{ bc } + \frac{ b }{ ac } + \frac{ c }{ ab }  \right) \\
                  &\mbox{}& \qquad \qquad \qquad - \frac{ 1 }{ 12 } ( a + b + c ) - \frac{ 1 }{ 12 } \left( \frac{ 1 }{ a } + \frac{ 1 }{ b } + \frac{ 1 }{ c }  \right)  \\
                  &\mbox{}& \qquad = \frac{ t^{ 2 }  }{ 2abc } - \frac{ t }{ 2 } \left( \frac{ 1 }{ ab } + \frac{ 1 }{ ac } + \frac{ 1 }{ bc }  \right) + \frac{ 1 }{ 12 } \left( \frac{ a }{ bc } + \frac{ b }{ ac } + \frac{ c }{ ab }  \right) \\
                  &\mbox{}& \qquad \qquad \qquad - \frac{ 1 }{ 12 } ( a + b + c ) + \frac{ 1 }{ 6 } \left( \frac{ 1 }{ a } + \frac{ 1 }{ b } + \frac{ 1 }{ c }  \right)  \ . \end{eqnarray*}
The larger zero of the right-hand side is an upper bound for the solution of the Frobenius problem: 
\begin{eqnarray*} &\mbox{}& f(a, b, c) \leq abc \left( \frac{ 1 }{ 2 } \left( \frac{ 1 }{ ab } + \frac{ 1 }{ bc } + \frac{ 1 }{ ac } \right) + \left[ \frac{ 1 }{ 4 } \left( \frac{ 1 }{ ab } + \frac{ 1 }{ bc } + \frac{ 1 }{ ac }  \right)^{ 2 } \right. \right. \\
                  &\mbox{}& \qquad \left. \left. - \frac{ 2 }{ abc } \left( \frac{ 1 }{ 12 } \left( \frac{ a }{ bc } + \frac{ b }{ ac } + \frac{ c }{ ab }  \right) - \frac{ 1 }{ 12 } ( a + b + c ) + \frac{ 1 }{ 6 } \left( \frac{ 1 }{ a } + \frac{ 1 }{ b } + \frac{ 1 }{ c }  \right)  \right)  \right]^{ 1/2 }    \right) \\
                  &\mbox{}& \quad \leq \frac{ 1 }{ 2 } \left( a + b + c \right) + abc \sqrt{ \frac{ 1 }{ 4 } \left( \frac{ 1 }{ ab } + \frac{ 1 }{ bc } + \frac{ 1 }{ ac }  \right)^{ 2 } + \frac{ 1 }{ 6 } \left( \frac{ 1 }{ ab } + \frac{ 1 }{ bc } + \frac{ 1 }{ ac } \right)  } \\
                  &\mbox{}& \quad = \frac{ 1 }{ 2 } \left( a + b + c \right) + abc \sqrt{ \frac{ 1 }{ 2 } \left( \frac{ 1 }{ ab } + \frac{ 1 }{ bc } + \frac{ 1 }{ ac } \right) \left( \frac{ 1 }{ 2 } \left( \frac{ 1 }{ ab } + \frac{ 1 }{ bc } + \frac{ 1 }{ ac } \right) + \frac{ 1 }{ 3 }  \right)  } \\
                  &\mbox{}& \quad \leq \frac{ 1 }{ 2 } \left( a + b + c \right) + abc \sqrt{ \frac{ 1 }{ 4 } \left( \frac{ 1 }{ ab } + \frac{ 1 }{ bc } + \frac{ 1 }{ ac } \right) } \ . \end{eqnarray*}
For the last inequality, we used the fact that $ \frac{ 1 }{ ab } + \frac{ 1 }{ bc } + \frac{ 1 }{ ac } \leq \frac{ 1 }{ 6 } + \frac{ 1 }{ 10 } + \frac{ 1 }{ 15 } = \frac{ 1 }{ 3 } $. 
This proves, using (\ref{threegen}) and (\ref{coprimegen}), 
\begin{theorem}\label{estimate} Let $ a_{1} \leq a_{2} \leq \dots \leq a_{n} $ be relatively prime. Then 
\[ f(a_{1}, \dots, a_{n}) \leq \frac{ 1 }{ 2 } \left( \sqrt{ a_{1} a_{2} a_{3} \left( a_{1} + a_{2} + a_{3} \right) } + a_{1} + a_{2} + a_{3} \right) + a_{4} + \dots + a_{n} \ .  \]
\end{theorem}
{\it Remarks.} 1. Sometimes the Frobenius problem is stated in a slightly different 
form: Given relatively prime positive integers $a_{1} , \dots , a_{n}$, 
find the largest value of $t$ such that $ \sum_{k=1}^{n} m_{k} a_{k} = t $ has no 
solution in \emph{nonnegative} integers $ m_{ 1 } , \dots , m_{ n } $. This number is denoted 
by $g(a_{1}, \dots, a_{n})$. It is, however, easy to see that 
  \[ g(a_{1}, \dots, a_{n}) = f(a_{1}, \dots, a_{n}) - a_{1} - \dots - a_{n} \ . \]
So we can restate Theorem \ref{estimate} in a more compact form as 
\[ g(a_{1}, \dots, a_{n}) \leq \frac{ 1 }{ 2 } \left( \sqrt{ a_{1} a_{2} a_{3} \left( a_{1} + a_{2} + a_{3} \right) } - a_{1} - a_{2} - a_{3} \right) \ .  \]

2. Bounds on the Frobenius number in the literature include results by 
Erd\H{o}s and Graham \cite{erdos} 
  \[ g(a_{1}, \dots, a_{n}) \leq 2 a_n \left\lfloor \frac {a_1} n \right\rfloor - a_1 \ , \] 
Selmer \cite{selmer} 
  \[ g(a_{1}, \dots, a_{n}) \leq 2 a_{n-1} \left\lfloor \frac {a_n} n \right\rfloor - a_n \ , \] 
and Vitek \cite{vitek} 
  \[ g(a_{1}, \dots, a_{n}) \leq \left\lfloor \frac 1 2 ( a_2 - 1 ) ( a_n - 2 ) \right\rfloor - 1 \ . \] 
Theorem \ref{estimate} is certainly of the same order. What might be 
more interesting, however, is the fact that the bound in Theorem \ref{estimate} is of 
a different nature than the bounds stated above: namely, it involves three variables, 
and is thus---especially in terms of estimating $g(a_1,a_2,a_3)$---more symmetric. 





\end{document}